\documentclass{amsart}
\usepackage{amsmath,amssymb}
\newtheorem{theorem}{Theorem}
\newtheorem{lemma}[theorem]{Lemma}
\newtheorem{proposition}[theorem]{Proposition}
\newtheorem{corollary}[theorem]{Corollary}

\begin{document}

\title[Strong cleanness of matrix rings]{Strong cleanness of matrix rings over commutative rings}
\author{Fran\c{c}ois Couchot}
\address{Laboratoire de Math\'ematiques Nicolas Oresme, CNRS UMR
  6139,
D\'epartement de math\'ematiques et m\'ecanique,
14032 Caen cedex, France}
\email{couchot@math.unicaen.fr}



\begin{abstract} Let $R$ be a commutative local ring. It is proved that  $R$ is Henselian if and only if each $R$-algebra which is a direct limit of module finite $R$-algebras is strongly clean. So, the matrix ring $\mathbb{M}_n(R)$ is strongly clean for each integer $n>0$ if $R$ is Henselian and we show that the converse holds if either the residue class field of $R$ is algebraically closed or $R$ is an integrally closed domain or $R$ is a valuation ring. It is also shown that each $R$-algebra which is locally a direct limit of module-finite algebras, is strongly clean if $R$ is a $\pi$-regular commutative ring.

\end{abstract}
\maketitle
As in \cite{Nic77} a ring $R$ is called \textbf{clean} if each element of $R$ is the sum of an idempotent and a unit. In \cite{HaNi01} Han and Nicholson proved that a ring $R$ is clean if and only if $\mathbb{M}_n(R)$ is clean for every integer $n\geq 1$. It is easy to check that each local ring is clean and consequently every matrix ring over a local ring is clean. On the other hand a ring $R$ is called \textbf{strongly clean} if each element of $R$ is the sum of an idempotent and a unit that commute. Recently, in \cite{WaCh04}, Chen and Wang gave an example of a commutative local ring $R$ with $\mathbb{M}_2(R)$ not strongly clean. This motivates the following interesting question: what are the commutative local rings $R$ for which $\mathbb{M}_n(R)$ is strongly clean for each integer $n\geq 1$? In \cite{CYZ06}, Chen, Yang and Zhou gave a complete characterization of commutative local rings $R$ with $\mathbb{M}_2(R)$ strongly clean. So, from their results and their examples, it is reasonable to conjecture that the Henselian rings are the only commutative local rings $R$ with $\mathbb{M}_n(R)$ strongly clean for each integer $n\geq 1$. In this note we give a partial answer to this problem. If $R$ is Henselian then $\mathbb{M}_n(R)$ is strongly clean for each integer $n\geq 1$ and the converse holds if $R$ is an integrally closed domain, a valuation ring or if its residue class field is algebraically closed.

All rings in this paper are associative with unity.
By \cite[Chapitre I]{Ray70} a commutative local ring $R$ is said to be \textbf{Henselian} if each commutative module-finite $R$-algebra is a finite product of local rings.  It was G. Azumaya (\cite{Azu51}) who first studied this property which was then developed by M. Nagata (\cite{Nag62}). The following theorem gives a new characterization of Henselian rings.

\begin{theorem}
\label{T:local} Let $R$ be a commutative local ring. Then the following conditions are equivalent:
\begin{enumerate}
\item $R$ is Henselian;
\item For each $R$-algebra $A$ which is a direct limit of module-finite algebras and for each integer $n\geq 1$, the matrix ring $\mathbb{M}_n(A)$ is strongly clean;
\item Each $R$-algebra $A$ which is a direct limit of module-finite algebras is clean.
\end{enumerate} 
\end{theorem}
\textbf{Proof.} $(1)\Rightarrow (2)$. Let $A$ be a direct limit of module-finite $R$-algebras and $a\in \mathbb{M}_n(A)$. Then $R[a]$ is a commutative module-finite $R$-algebra. Since $R$ is Henselian, $R[a]$ is a finite direct product of local rings. So $R[a]$ is  clean. Hence $a$ is a sum of an idempotent and a unit that commute.

It is obvious that $(2)\Rightarrow (3)$.

$(3)\Rightarrow (1)$. Let $A$ be a commutative module-finite $R$-algebra and let $J(A)$ be its Jacobson radical. Since $J(R)A\subseteq J(A)$, where $J(R)$ is the Jacobson radical of $R$, we deduce that $A/J(A)$ is semisimple artinian. By \cite[ Propositions 1.8 and 1.5]{Nic77} idempotents can be lifted modulo $J(A)$.  Hence $A$ is semi-perfect. It follows that $A$ is a finite product of local rings, whence $R$ is Henselian.
\qed

\medskip Let $\mathcal{P}$ be a ring property. We say that an algebra $A$ over a commutative ring $R$ is  \textbf{locally $\mathcal{P}$} if $A_P$ satisfies $\mathcal{P}$ for each maximal ideal $P$ of $R$.

\begin{corollary} \label{C:global}
Let $R$ be a  commutative ring. Then the following conditions are equivalent:
\begin{enumerate}
\item $R$ is clean and locally Henselian;
\item For each $R$-algebra $A$ which is locally a direct limit of module-finite algebras and for each integer $n\geq 1$,  $\mathbb{M}_n(A)$ is strongly clean;
\item Each $R$-algebra $A$ which is  locally a direct limit of module-finite algebras is clean.
\end{enumerate} 
\end{corollary}
\textbf{Proof.} $(1)\Rightarrow (2)$. Let $A$ be an $R$-algebra which is locally a direct limit of module-finite algebras and $a\in\mathbb{M}_n(A)$. Consider the following polynomial equations: $E+U=a,\ E^2=E,\ UV=1,\ VU=1,\ EU=UE$. By Theorem~\ref{T:local} these equations have a solution in $\mathbb{M}_n(A_P)$, for each maximal ideal $P$ of $R$. So, by \cite[Theorem I.1]{Cou07} they have a solution in $\mathbb{M}_n(A)$ too. 

It is obvious that $(2)\Rightarrow (3)$.

$(3)\Rightarrow (1)$. Let $P$ be a maximal ideal of $R$ and let $A$ be a module-finite $R_P$-algebra. Since $R$ is clean, the natural map $R\rightarrow R_P$ is surjective by \cite[Theorem I.1 and Proposition III.1]{Cou07}. So $A$ is a module-finite $R$-algebra. It follows that $A$ is  clean. By Theorem~\ref{T:local} $R_P$ is Henselian.
\qed

\medskip 
A ring $R$ is said to be \textbf{strongly $\pi$-regular} if, for each $r\in R$, there exist $s\in R$ and an integer $q\geq 1$ such that $r^q=r^{q+1}s$.

\begin{corollary}
\label{C:pireg} Let $R$ be a strongly $\pi$-regular commutative ring. Then, for each $R$-algebra $A$ which is locally a direct limit of module-finite algebras and for each integer $n\geq 1$, the matrix ring $\mathbb{M}_n(A)$ is strongly clean.
\end{corollary}
\textbf{Proof.} It is known that $R$ is clean and that each prime ideal is maximal. So, for every maximal $P$, $PR_P$ is a nilideal of $R_P$. Hence $R_P$ is Henselian. We conclude by Corollary~\ref{C:global}. \qed

\medskip  By \cite[Th\'eor\`eme 1]{Dis76} each strongly $\pi$-regular $R$ satisfies the following condition: for each $r\in R$, there exist $s\in R$ and an integer $q\geq 1$ such that $r^q=sr^{q+1}$. Moreover, by \cite[Proposition 2.6.iii)]{BuMe88} each strongly $\pi$-regular ring is strongly clean. So, Corollary~\ref{C:pireg} is also a consequence of the following proposition. (Probably, this proposition is already known).

\begin{proposition} \label{P:pireg}
Let $R$ be a strongly $\pi$-regular commutative ring. Then, for each $R$-algebra $A$ which is locally a direct limit of module-finite algebras and for each integer $n\geq 1$, the matrix ring $\mathbb{M}_n(A)$ is strongly $\pi$-regular.
\end{proposition}
\textbf{Proof.} Let $S=\mathbb{M}_n(A)$ and $s\in S$. Then $R[s]$ is locally a module-finite algebra. It is easy to prove that each prime ideal of $R[s]$ is maximal. Consequently $R[s]$ is strongly $\pi$-regular. So, $S$ is strongly $\pi$-regular too. \qed

\medskip 
The following lemma will be useful in the sequel.

\begin{lemma}
\label{L:polyReduc} Let $R$ be a commutative local ring with maximal ideal $P$. Let $n$ be an integer $>1$ such that $\mathbb{M}_n(R)$ is strongly clean. Let $f$ be a monic polynomial of degree $n$ with coefficients in $R$ such that $f(0)\in P$ and $f(a)\in P$ for some $a\in R\setminus P$. Then $f$ is reducible.
\end{lemma}
\textbf{Proof}. Let $A\in\mathbb{M}_n(R)$ such that its characteristic polynomial is $f$, i.e. $f=\det(XI_n-A)$, where $I_n$ is the unit element of $\mathbb{M}_n(R)$. Then $A=E+U$ where $E$ is idempotent, $U$ is invertible and $EU=UE$. First we assume that $a=1$. So, $0$ and $1$ are eigenvalues of $\overline{A}$ the reduction of $A$ modulo $P$. Consequently $A$ and $A-I_n$ are not invertible.  It follows that $E\ne I_n$ and $E\ne 0_{n,n}$ where $0_{p,q}$ is the $p\times q$ matrix whose coefficients are $0$. Let $F$ be a free $R$-module  of rank $n$ and let $\epsilon$ be the endomorphism  of $F$ for which $E$ is the matrix associated  with respect to some basis. Then $F=\mathrm{Im}\ \epsilon\oplus\mathrm{Ker}\ \epsilon$. Moreover $\mathrm{Im}\ \epsilon$ and $\mathrm{Ker}\ \epsilon$ are free because $R$ is local. Consequently there exists a $n\times n$ invertible matrix $Q$ such that:
\[QEQ^{-1}=B=\begin{pmatrix}
I_p & 0_{p,q} \\
0_{q,p}& 0_{q,q}
\end{pmatrix}
\]
where $p$ is an integer such that $0<p<n$ and $q=n-p$. Since $E$ and $A$ commute, then $B$ and $QAQ^{-1}$ commute too. So, $QAQ^{-1}$ is of the form:
\[QAQ^{-1}=\begin{pmatrix}
C & 0_{p,q} \\
0_{q,p}& D
\end{pmatrix}
\]
where $C$ is a $p\times p$ matrix and $D$ is a $q\times q$ matrix. We deduce that $f$ is the product of the characteristic polynomial $g$ of $C$ with the characteristic polynomial $h$ of $D$. Let us observe that $(C-I_p)$ and $D$ are invertible. So, $g(1)\notin P,\ h(0)\notin P,\ g(0)\in P$ and $h(1)\in P$. Now suppose that $a\ne 1$.  Then $a^{-n}f(X)=g(Y)$ where $Y=a^{-1}X$ and $g$ is a monic polynomial of degree $n$. We easily check that $g(1)\in P$ and $g(0)\in P$. It follows that $g$ is reducible, whence $f$ is reducible too. \qed

\medskip A commutative ring $R$ is a \textbf{valuation ring} (respectively \textbf{arithmetic}) if its lattice of ideals is totally ordered by inclusion (respectively distributive).
\begin{theorem}
\label{T:matrix} Let $R$ be a local commutative ring with maximal ideal $P$ and with residue class field $k$. Consider the following two conditions:
\begin{enumerate}
\item $R$ is Henselian;
\item  the matrix ring $\mathbb{M}_n(R)$ is strongly clean $\forall n\in\mathbb{N}^*$.
\end{enumerate} 
Then $(1)\Rightarrow (2)$ and the converse holds if $R$ satisfies one of the following properties:
\begin{itemize}
\item[(a)] $k$ is algebraically closed;
\item[(b)] $R$ is an integrally closed domain;
\item[(c)] $R$ is a valuation ring.
\end{itemize} 
\end{theorem}
\textbf{Proof.} By Theorem~\ref{T:local} it remains to prove that $(2)$ implies $(1)$ when one of (a), (b) or (c) is valid. We will use \cite[Theorem 1.4]{Bon02} and \cite[Theorem II.7.3.(iv)]{FuSa01}. Consider the  polynomial $f=X^n+c_{n-1}X^{n-1}+\dots +c_1X+c_0$ and assume that $\exists m,\ 1\leq m<n$ such that $c_m\notin P$ and $c_i\in P,\ \forall i<m$. Since $c_0\in P$, we see that $f(0)\in P$.

Hence, if $k$ is algebraically closed, $\exists a\in R\setminus P$ such that $f(a)\in P$. By Lemma~\ref{L:polyReduc} $f$ is reducible. So, by \cite[Theorem 1.4]{Bon02} $R$ is Henselian. 

If $R$ is an integrally closed domain, we take $m=n-1$ for proving the condition (iv) of \cite[Theorem II.7.3]{FuSa01}. In this case $f(-c_{n-1})\in P$. By Lemma~\ref{L:polyReduc} (possibly applied several times) $f$ satisfies the condition (iv) of \cite[Theorem II.7.3]{FuSa01}. Hence $R$ is Henselian. 

Assume that $R$ is a valuation ring. Let $N$ be the nilradical of $R$ and let $R'=R/N$. We know that $R$ is Henselian if and only if $R'$ is Henselian too. For each $n\in\mathbb{N}^*$, $\mathbb{M}_n(R')$ is strongly clean. Since $R'$ is a valuation domain, $R'$ is integrally closed. It follows that $R'$ and $R$ are Henselian. \qed

\begin{corollary} \label{C:arithmetic}
Let $R$ be an arithmetic commutative ring. Then the following conditions are equivalent:
\begin{enumerate}
\item $R$ is clean and locally Henselian;
\item the matrix ring $\mathbb{M}_n(R)$ is strongly clean $\forall n\in\mathbb{N}^*$.
\end{enumerate} 
\end{corollary}
\textbf{Proof.} By Corollary~\ref{C:global} it remains to show $(2)\Rightarrow (1)$. Let $P$ be a maximal ideal of $R$. Since $R$ is clean the natural map $R\rightarrow R_P$ is surjective by \cite[Theorem I.1 and Proposition III.1]{Cou07}. So, $\mathbb{M}_n(R_P)$ is strongly clean $\forall n\in\mathbb{N}^*$. Theorem~\ref{T:matrix} can be applied  because $R_P$ is a valuation ring. We conclude that $R_P$ is Henselian. \qed

\medskip  

The following generalization of \cite[Theorem 8]{CYZ06} holds even if the properties $(a),(b),(c)$ of Theorem~\ref{T:matrix} are not satisfied.

\begin{theorem} \label{T:n5}
Let $R$ be a local commutative ring with maximal ideal $P$ and with residue class field $k$. Let $p$ be an integer such that  $2\leq p\leq 5$. Then the following conditions are equivalent:
\begin{enumerate}
\item $\mathbb{M}_n(R)$ is strongly clean $\forall n,\ 2\leq n\leq p$;
\item  each  monic polynomial $f$ of degree $n,\ 2\leq n\leq p$, for which $f(0)\in P$ and $f(1)\in P$, is reducible. 
\end{enumerate}  
\end{theorem}
\textbf{Proof.} By Lemma~\ref{L:polyReduc} it remains to prove that $(2)\Rightarrow (1)$. Let $A\in\mathbb{M}_n(R)$. We denote by $f$ the characteristic polynomial of $A$. If $A$ is invertible then $A=0_{n,n}+A$. If $A-I_n$ is invertible then $A=I_n+(A-I_n)$. So, we may assume that $A$ and $(A-I_n)$ are not invertible. It follows that $f(0)\in P$ and $f(1)\in P$. Then, $f=gh$ where $g$ and $h$ are monic polynomials of degree $\geq 1$. We may assume that $g(0)\in P,\ g(1)\notin P,\ h(0)\notin P$ and $h(1)\in P$ (possibly by applying condition $(2)$ several times). We denote by $\bar{f},\ \bar{g},\ \bar{h}$ the images of $f,\ g,\ h$ by the natural map $R[X]\rightarrow k[X]$. If $\bar{g}$ and $\bar{h}$ have a common factor of degree $\geq 1$ then this factor is of degree $1$ because $n\leq 5$. In this case $\exists a\in R\setminus P$ such that $g(a)\in P$ and $h(a)\in P$. As in the proof of Lemma~\ref{L:polyReduc} we show that $g$ is reducible. Hence, after changing $g$ and $h$, we get that $\bar{g}$ and $\bar{h}$ have no common divisor of degree $\geq 1$. It follows that there exist two polynomials $u$ and $v$  with coefficients in $R$ such that $\bar{u}\bar{g}+\bar{v}\bar{h}=1$. Since $PR[A]$ is contained in the Jacobson radical of $R[A]$, we may assume that $u(A)g(A)+v(A)h(A)=I_n$. We put $e=vh$. Then we easily check that $e(A)$ is idempotent. It remains to show that $(A-e(A))$ is invertible. It is enough to prove that $(\bar{A}-\bar{e}(\bar{A}))$ is invertible because $P\mathbb{M}_n(R)$ is the Jacobson radical of $\mathbb{M}_n(R)$. Let $V$ be a vector space of dimension $n$ over $k$ and let $\mathcal{B}$ be a basis of $V$. Let $\alpha$ be the endomorphism of $V$ for which  $\bar{A}$ is the matrix associated with respect to $\mathcal{B}$. We put $\epsilon=\bar{e}(\alpha)$. Since $V$ has finite dimension, it is sufficient to show that $(\alpha-\epsilon)$ is injective. Let $w\in V$ such that $\alpha(w)=\epsilon(w)$. It follows that $\alpha(\epsilon(w))=\epsilon(\alpha(w))=\epsilon^2(w)=\epsilon(w)$. Since $\bar{e}$ is divisible by $(X-\bar{1})$ we get that $\epsilon(w)=0$. So, $\alpha(w)=0$. We deduce that $\epsilon(w)=w$ because $\bar{e}-\bar{1}$ is divisible by $X$. Hence $w=0$. \qed

\end{document}